\documentclass[12pt]{article}

\usepackage[centertags]{amsmath}
\usepackage{amsfonts}
\usepackage{amssymb}

\begin{document}

\newcounter{abcd}
\setcounter{abcd}{0}\Alph{abcd}

\newtheorem{pred0}{Theorem}
\newtheorem{pred01}{Theorem}
\newtheorem{predA}{Theorem}
\renewcommand{\thepredA}{\Alph{predA}}

\newtheorem{propo}{Proposition}

\newtheorem{predB}{Theorem B}

\newcommand {\Hol}{\mathop{\rm Hol}\nolimits}

\newcommand{\G}{\mathcal{G}}
\newtheorem{predD}{Definition}
\newtheorem{prop}{Proposition}
\newtheorem{lemma}{Lemma}
\newtheorem{lemma1}{Lemma}
\newtheorem{corol}{Corollary}
\newcommand{\pr}{\noindent{\bf Proof.}\quad }
\newcommand{\prr}{\noindent{\bf Proof of Theorem 1.}\quad }
\newcommand{\prrr}{\noindent{\bf Proof of Theorem 3.}\quad }
\newcommand{\epr}{\ $\blacksquare$  \vspace{2mm} }
\newcommand{\eprr}{$\blacksquare$ }
\newtheorem{rem}{Remark}
\newcommand{\prone}{\noindent{\bf Proof of Theorem 1.}\quad }
\newcommand{\ac}{\noindent{\bf Acknowledgment.}\quad }
\renewcommand{\Re}{\mathop{\rm Re}\nolimits}
\renewcommand{\Im}{\mathop{\rm Im}\nolimits}

\title{Rigidity of holomorphic generators and one-parameter semigroups}

\author{Mark Elin
\\ {\small Department of Mathematics, ORT  Braude College,}
\\ {\small P.O. Box 78, Karmiel 21982, Israel}
\\ {\small e-mail: mark.elin@gmail.com}
\\ Marina Levenshtein
\\ {\small Department of Mathematics,}
\\ {\small The Technion --- Israel Institute of Technology}
\\ {\small 32000 Haifa, Israel}
\\ {\small e-mail: marlev@list.ru}
\\ David Shoikhet
\\ {\small Department of Mathematics, ORT Braude College,}
\\ {\small P.O. Box 78, Karmiel 21982, Israel}
\\ Roberto Tauraso
\\ {\small Dipartimento di Matematica, Universit\`{a} di Roma `Tor Vergata',} \\
{\small Via della Ricerca Scientifica, 00133 Roma, Italy}
\\ {\small e-mail: tauraso@mat.uniroma2.it}
}
\date{ }
\maketitle

\begin{abstract}
In this paper we establish a rigidity property of holomorphic
generators by using their local behavior at a boundary point
$\tau$ of the open unit disk $\Delta$. Namely, if
$f\in\mathrm{Hol}(\Delta,\mathbb{C})$ is the generator of a
one-parameter continuous semigroup $\{F_{t}\}_{t\geq0}$, we state
that the equality $f(z)=o\left(|z-\tau|^{3}\right)$ when
$z\rightarrow\tau$ in each non-tangential approach region at
$\tau$ implies that $f$ vanishes identically on $\Delta$. Note,
that if $F$ is a self-mapping of $\Delta$ then $f=I-F$ is a
generator, so our result extends the boundary version of the
Schwarz Lemma obtained by D. Burns and S. Krantz. We also prove
that two semigroups $\{F_{t}\}_{t\geq0}$ and $\{G_{t}\}_{t\geq0}$,
with generators $f$ and $g$ respectively, commute if and only if
the equality $f=\alpha g$ holds for some complex constant
$\alpha$. This fact gives simple conditions on the generators of
two commuting semigroups at their common null point $\tau$ under
which the semigroups coincide identically on $\Delta$.
\end{abstract}

\section{Introduction.}

Let $\Delta= \{ z \in \mathbb{C}:\ |z|<1 \} $ be the open unit
disk in the complex plane $\mathbb{C}$, and let $H= \{ z \in
\mathbb{C}:\ \Re z> 0 \} $ be the right half-plane. We denote by
$\mathrm{Hol}(\Delta,D)$ the set of all holomorphic functions on
$\Delta$ which map $\Delta$ into a set $D\subset \mathbb{C}$, and
by $\mathrm{Hol}(\Delta)$ the set of all holomorphic self-mappings
of $\Delta$, i.e.,
$\mathrm{Hol}(\Delta)=\mathrm{Hol}(\Delta,\Delta)$.

The problem of finding conditions for a holomorphic function $F$
to coincide identically with a given holomorphic function $G$ when
they have a similar behavior on some subset of
$\overline{\Delta}$, has been studied by many mathematicians.

The following assertions are classical:

$\bullet$ \emph { If $F$ and $G$ are holomorphic in $\Delta$ and
$F=G$ on a subset of $\Delta$ that has a nonisolated point, then
$F\equiv G$ on $\Delta$ (Vitali's uniqueness principle).}

$\bullet$ \emph { If $F$ and $G$ are holomorphic in $\Delta$ and
continuous on $\overline{\Delta}$, and $F=G$ on some arc $\gamma$
of the boundary $\partial\Delta$, then $F\equiv G$ on $\Delta$.}

In the point of view of complex dynamics it is natural to study
conditions on derivatives of $F$ and $G$ at specific points to
conclude that $F\equiv G$.

If, for example, $G$ is the identity mapping $I$ and
$\tau\in\Delta$ is the Denjoy--Wolff point of
${F\in\Hol(\Delta)}$, then the equalities $F(\tau)=G(\tau)$ and
$F'(\tau)=G'(\tau)$ provide $F\equiv G$ by the Schwarz Lemma. The
same conclusion holds for an arbitrary holomorphic function $G$ on
$\Delta$, if $F$ commutes with $G$ and satisfies the conditions
$F(\tau)=G(\tau)=\tau$ and $F'(\tau)=G'(\tau)\neq 0$ (see, for
instance, \cite{C-M}, \cite{B-T-V}).

Different ``identity principles" have been recently studied by
several mathematicians under suitable boundary conditions. In
general, the following three cases are considered.

\quad  (A) \quad    $G$ is the identity mapping;

\quad  (B) \quad    $G$ is an arbitrary self-mapping of $\Delta$,
and $F$  commutes with $G$,\linebreak i.e., $F \circ G =G \circ
F$;

\quad  (C) \quad    $G$ is a constant mapping.

Regarding Case A the following result is due to D.~Burns and
S.~Krantz.

\begin{predA}[\cite{B-K}]
Let $F\in\mathrm{Hol}(\Delta)$ and
\begin{equation}\label{a1}
 F(z)=1+(z-1)+O\left((z-1)^{4}\right).
\end{equation}
Then $F\equiv I$.
\end{predA}

For Case B a uniqueness theorem was given by R. Tauraso in
\cite{T} (see also \cite{B-T-V}). To formulate this result we need
the following notation. Let $F\in\Hol(\Delta)$ and
$\tau\in\partial\Delta$. We say that $F\in C^m_K(\tau)$ if it
admits the following representation
\[
F(z)=\tau +F'(\tau)(z-\tau)+\ldots
+\frac{F^{(m)}(\tau)}{m!}(z-\tau)^m +o\left(|z-\tau|^m\right)
\]
when $z\to\tau$ in each non-tangential approach region at $\tau$.
Moreover, we say that $F\in C^m(\tau)$ if the limit is taken in
the full disk.

\begin{predA}[\cite{T}]
Let $F,G \in\mathrm{Hol}(\Delta)$ be commuting functions with a
common Denjoy--Wolff point $\tau\in\partial\Delta$. If one of the
following conditions holds then $F\equiv G$.

(i) $F'(\tau)=G'(\tau)<1$;

(ii) $F\in C^{2}(\tau),\ G\in C^2_K(\tau),\ F'(\tau)=1,\
F''(\tau)=G''(\tau)\neq 0$ and $\Re \tau F''(\tau)>0$;

(iii) $F,\, G\in C^{2}(\tau),\ F'(\tau)=1,\
F''(\tau)=G''(\tau)\neq 0$ and $\Re \tau F''(\tau)=0$;

(iv) $F\in C^{3}(\tau),\ G\in C^{3}_{K}(\tau),\ F'(\tau)=1,\
F''(\tau)=G''(\tau)=0$ and $F'''(\tau)=G'''(\tau)$.
\end{predA}

For Case C, when $G$ is a constant mapping, the following fact is
an immediate consequence of the Julia--Wolff--Carath\'eodory
Theorem.

$\bullet$ \emph { If $F \in\Hol(\Delta,\overline{\Delta})$, then
the conditions $\lim\limits_{r\rightarrow 1^{-}} F(r \tau)=\tau$
and ${\lim\limits_{r\rightarrow 1^{-}} F'(r \tau)=0}$ for some
$\tau\in\partial\Delta$ imply that $F\equiv\tau$.}

In fact, the considering of holomorphic functions $f$ which are
not necessarily self-mappings  is more relevant in this situation.
Various results in this direction were established by
S.~Migliorini and F.~Vlacci in~\cite{M-V}.

In what follows we denote by symbol $\angle\lim\limits_{z\to\tau}$
the angular limit of a function defined in $\Delta$ at a boundary
point $\tau\in\partial\Delta$.

\begin{predA}[see \cite{M-V}] Let $\tau\in\partial\Delta$.

If $f \in\mathrm{Hol}(\Delta,\overline{H})$, then
\begin{equation}\label{a4}
\angle\lim\limits_{z\rightarrow \tau}\frac{f(z)}{z-\tau}=0
\end{equation}
implies that $f\equiv0$.

More general, if $f \in\mathrm{Hol}(\Delta,\mathbb{C})$, and
$f(\Delta)$ is contained in a wedge of angle $\pi\alpha$, $
0<\alpha \leq 2$, with vertex at the origin, then the condition
\begin{equation}\label{a6}
\angle\lim\limits_{z\rightarrow \tau}\frac{f(z)}{(z-\tau)
^{\alpha}}=0
\end{equation}
implies that $f\equiv0$.
\end{predA}

Although the classes $\Hol(\Delta)$ of holomorphic self-mappings
of $\Delta$ and $\Hol(\Delta,H)$ of functions with positive real
part are connected by the composition with the Cayley transform,
Theorem~A is not a direct consequence of Theorem~C, and
conversely.

In this note we find rigidity principles for some classes of
holomorphic functions produced by continuous dynamical systems,
which are related to both $\Hol(\Delta)$ and $\Hol(\Delta,H)$. In
particular, by this way one can establish a bridge between
Theorems A and C.

We consider, inter alia, the class of mappings
$F\in\mathrm{Hol}(\Delta,\mathbb{C})$ which are continuous on
$\overline{\Delta}$ and satisfy the boundary flow-invariance
condition
\begin{equation}\label{a7} \Re F(z)\overline{z}\leq 1, \, z\in \partial \Delta .
\end{equation}
In particular, each function $F\in \mathrm{Hol}(\Delta)$ which is
continuous on $\overline{\Delta}$ belongs to this class.

Condition (\ref{a7}) can be rewritten in the form
\begin{equation}\label{a8}
\Re f(z)\overline{z}\geq 0, \, z\in \partial \Delta ,
\end{equation}
where
\begin{equation}\label{a9} f(z)=z-F(z).
\end{equation}

Note that each mapping $f$ satisfying (\ref{a8}) belongs to the
class $\G(\Delta)$ of so-called semigroup generators on $\Delta$.

Our main purpose is to establish boundary conditions for a
function ${f\in\G(\Delta)}$ to vanish on $\Delta$ identically.

First, we recall that a family $S=\left\{F_t\right\}_{t\geq
0}\subset\mathrm{Hol}(\Delta)$ is said to be {\bf a one-parameter
continuous semigroup on $\Delta$} if

(i) $F_{t}(F_{s}(z))=F_{t+s}(z)$ for all $t,s\geq 0,$

(ii) $\lim\limits_{t\rightarrow 0^+}F_{t}(z)=z$ for all $z\in
\Delta$.

Furthermore, it follows by a result of E. Berkson and H. Porta
\cite{B-P} that each continuous semigroup is differentiable in
$t\in\mathbb{R}^{+}=[0,\infty ),$ (see also \cite{AM-88} and
\cite{R-S-98}). So, for each continuous semigroup
$S=\left\{F_{t}\right\}_{t\geq 0}\subset\mathrm{Hol}(\Delta)$, the
limit
\begin{equation}
\lim_{t\rightarrow 0^{+}}\frac{z-F_{t}(z)}{t}=f(z),\quad
z\in\Delta, \label{1a}
\end{equation}
exists and defines a holomorphic mapping
$f\in\mathrm{Hol}(\Delta,\mathbb{C})$. This mapping $f$ is called
the {\bf (infinitesimal) generator of}
$S=\left\{F_{t}\right\}_{t\geq 0}.$ Moreover, the function
$u(=u(t,z)),\ (t,z)\in\mathbb{R}^{+}\times\Delta$, defined by
$u(t,z)=F_{t}(z)$ is the unique solution of the Cauchy problem
\begin{equation}\label{2a}
\left\{
\begin{array}{l}
{\displaystyle\frac{\partial u(t,z)}{\partial t}}+f(u(t,z))=0,\vspace{3mm}\\
u(0,z)=z,\quad z\in \Delta.
\end{array}
\right.
\end{equation}
The class of all holomorphic generators on $\Delta$ is denoted by
$\G(\Delta)$.

Note, that if $F\in\mathrm{Hol}(\Delta)$, then the function
$f=I-F$ belongs to $\G(\Delta)$ (see Corollary 3.3.1 in \cite{S}).

The following assertion combines characterizations of the class
$\G (\Delta)$ obtained in \cite{A-E-R-S}, \cite{A-R-S}  and
\cite{B-P}.

\begin{prop}
Let  $f \in \mathrm{Hol}(\Delta,\mathbb{C})$. The following are
equivalent:

(i) $f$ is a semigroup generator on $\Delta$;

(ii) $\Re f(z) \overline{z} \geq \Re f(0) \overline{z} \left
(1-|z|^{2} \right) \; $ for all  $\; z\in\Delta$;

(iii) there exists a unique point $\tau\in\overline{\Delta}$ such
that
\begin{equation}\label{a43}
f(z)=(z-\tau )(1-\overline{\tau}z)g(z),\;\; z\in\Delta,
\end{equation}
where $g\in\mathrm{Hol}(\Delta,\mathbb{C})$, $\Re g(z)\geq 0$.

(iv) $f$ admits the representation $$f(z)=a-\overline{a}z^{2}+z
p(z),$$ where $a\in\mathbb{C}$ and
$p\in\mathrm{Hol}(\Delta,\mathbb{C})$ with $\Re p(z)\geq 0$.
\end{prop}

\begin{rem}
The point $\tau$ in (\ref{a43}) is the Denjoy--Wolff point of the
semigroup $\left \{ F_{t} \right \}_{t\geq 0} $ generated by $f$.
If $\tau\in\Delta$ then $f(0)=0$ and $\Re f'(\tau)\geq 0$. If
$\tau\in\partial\Delta$ then the angular limit
$\angle\lim\limits_{z\rightarrow
\tau}\frac{f(z)}{z-\tau}=:f'(\tau)$ exists and is a nonnegative
real number (see \cite{E-S}).
\end{rem}

\section{Rigidity of infinitesimal generators.}

\begin{pred0}
Let $f\in \G(\Delta)$. Suppose that for some
$\tau\in\partial\Delta$
\[
f(z)=a(z-\tau)^{3}+o\left(|z-\tau|^{3}\right)
\]
when $z\to\tau$ in each  non-tangential approach region at $\tau$.
Then $a\tau^{2}$ is a nonnegative real number. Moreover, $a=0$ if
and only if $f\equiv 0$.
\end{pred0}
To prove Theorem 1 we need the following lemma.

\begin{lemma}
Let $g\in\mathrm{Hol}(\Delta,\overline H)$. Then for each
$\tau\in\partial\Delta$ the limit
\begin{equation}\label{b1}
k=\angle\lim\limits_{z\rightarrow \tau}
\frac{g(z)}{1-\overline{\tau} z}
\end{equation}
is either a nonnegative real number or infinity. Moreover,
$g\equiv 0$ if and only if $k=0$.
\end{lemma}

\pr Denote by $C_{\tau}(z)=\frac{\tau-z}{\tau+z}$ the Cayley
transform and set $h=C_{\tau}^{-1}\circ
g\in\mathrm{Hol}(\Delta,\overline{\Delta})$. By the
Julia--Wolff--Carath\'eodory theorem the limit
$$\beta_{h}=\angle\lim\limits_{z\rightarrow \tau}
\frac{\tau-h(z)}{\tau-z}$$ exists and is either a nonnegative real
number or infinity. Moreover, $\beta_{h}=0$ if and only if
$h\equiv \tau$.

For any $z\in\Delta$ we have
\begin{equation}\label{b2}
\frac{g(z)}{1-\overline{\tau} z
}=\frac{\tau-h(z)}{\tau-z}\,\cdot\, \frac{\tau}{\tau + h(z)}\,.
\end{equation}
Hence, $k=0$ if and only if $\beta_{h}=0$, and therefore $g\equiv
0$.

If $\beta_{h}$ is a positive real number, $\beta_{h}>0$,  then
$\angle\lim\limits_{z\rightarrow \tau}h(z)=\tau$ and,
consequently, $$k=\angle\lim\limits_{z\rightarrow \tau}
\frac{\tau-h(z)}{\tau-z}\,\cdot\,\angle\lim\limits_{z\rightarrow
\tau}\frac{\tau}{\tau+h(z)} =\frac{\beta_{h}}{2}>0.$$

Let $\beta_{h}=\infty$. Since $\Re
\frac{\tau}{\tau+h(z)}\geq\frac{1}{2}\,$, formula (\ref{b2})
implies that $k=\infty$. \epr

\noindent{\bf Alternative proof.} If $g\neq 0$, then the function
$p$ defined by $p(z):=\frac{1}{g(z)}$ belongs to
$\mathrm{Hol}(\Delta,H)$. It is easy to see that for all
$\zeta\in\partial\Delta$ the expression
$\frac{(1-z\overline{\tau})(1+z\overline{\zeta})}{1-z\overline{\zeta}}$
is bounded on each  non-tangential approach region at $\tau$. Then
it follows by the Riesz--Herglots formula that
\[
\angle\lim\limits_{z\rightarrow \tau} (1-z\bar\tau)p(z)=
\angle\lim\limits_{z\rightarrow \tau}\oint
\limits_{\partial\Delta}\frac{(1-z\overline{\tau})(1+z\overline{\zeta})}{1-z\overline{\zeta}}dm_{p}(\zeta)
=2m_{p}(\tau)\ge0,
\]
where $dm_{p}$ is a probability measure on $\partial\Delta$.
Setting $k=\frac{1}{2m_{p}(\tau)}$ we get our assertion. \epr

\prr Since $$\angle\lim\limits_{z\rightarrow
\tau}\frac{f(z)}{z-\tau}=0, $$ it follows from \cite{E-S} that
$\tau\in\partial\Delta$ is the Denjoy--Wolff point for the
semigroup $\left\{F_{t}\right\}_{t\geq 0 }$ generated by $f$. Then
by Proposition 1 the function $f$ admits the representation
(\ref{a43}): $$f(z)=(z-\tau)(1-z\overline{\tau})g(z)$$ with some
$g\in\mathrm{Hol}(\Delta,\overline{H})$. Hence, by Lemma 1 $$a
\tau^{2}=\tau^{2}\angle\lim\limits_{z\rightarrow
\tau}\frac{f(z)}{(z-\tau)^{3}}= \angle\lim\limits_{z\rightarrow
\tau}\frac{g(z)}{1-\overline{\tau} z}=k\geq 0.$$ Obviously, $a=0$
if and only if $k=0$. In this case $g\equiv 0$, so $f\equiv 0$.
\eprr

\begin{corol}[cf. Theorem 5 in \cite{B-T-V}.]
Let $F\in\mathrm{Hol}(\Delta,\mathbb{C})$ be continuous on
$\overline{\Delta}$  and satisfy the boundary condition
\[
\Re F(z)\overline{z}\leq 1, \, z\in\partial\Delta.
\]
If $F$ admit the representation
\[
F(z)=\tau+(z-\tau)+b(z-\tau)^{3}+o\left(|z-\tau|^{3}\right)
\]
when $z\to\tau$ in each  non-tangential approach region at some
point $\tau\in\partial\Delta$, then $b \tau^{2}\leq 0$. Moreover,
$b=0$ if and only if $F\equiv I$.
\end{corol}

As a consequence of Lemma 1 we also obtain the following
assertion.
\begin{corol}
Let  $f \in \G (\Delta)$ be such that $f(\tau)=0$ for some
$\tau\in\partial\Delta$ and $f(0)=a\in\mathbb{C}$. Suppose that
$f$ has a finite angular derivative at $\tau$. Then $f'(\tau)$ is
a real number with $f'(\tau)\leq -2\Re (\overline{a}\tau)$.
Moreover, $f'(\tau)=- 2\Re (\overline{a}\tau) \;$ if and only if
$f$ generates a group of automorphisms.
\end{corol}

\pr
By Proposition 1 (iv) $f$ admits the representation
\begin{equation}\label{a101}
f(z)=a-\overline{a}z^{2}+z p(z), \, \,  z\in\Delta,
\end{equation}
where $p\in \mathrm{Hol}(\Delta,\mathbb{C})$ with $\Re p(z)\geq
0$.

Since $f(\tau)=0$, we have
$p(\tau)=\overline{a}\tau-a\overline{\tau}=2i\Im
(\overline{a}\tau) \,$ is pure imaginary.

Then it follows from (\ref{a101}), that
$$f'(\tau)=\angle\lim\limits_{z\rightarrow
\tau}\frac{a-\overline{a}z^{2}+z p(z)}{z-\tau}=-2\Re
(a\overline{\tau})+\angle\lim\limits_{z\rightarrow
\tau}\frac{p(z)-2i\Im (\overline{a}\tau)}{z\overline{\tau}-1}\,.$$

Applying Lemma 1 to the function $g(z)=p(z)-2i\Im
(\overline{a}\tau)$, we get \linebreak $f'(\tau)\leq -2\Re(
a\overline{\tau})$.

Moreover, $f'(\tau)=-2\Re( a\overline{\tau}) \,$ if and only if
$p\equiv 2i\Im (\overline{a}\tau)$, i.e., $f(z)=a+2i\Im
(\overline{a}\tau)\cdot z-\overline{a}z^{2} \,$.

By Proposition 3.5.1 in \cite{S} each function of the form
$f(z)=a+ib z-\overline{a}z^{2}$, with $a\in\mathbb{C}$ and
$b\in\mathbb{R}$, generates a group of automorphisms of $\Delta$.
The proof is complete.
\epr

\begin{corol}
Let $F\in\mathrm{Hol}(\Delta)$ be such that $F(\tau)=\tau$ and
$F(0)=a, \, a\in\Delta$. Suppose that $F$ has a finite angular
derivative at $\tau$. Then $F'(\tau)\geq 1+2\Re
(\overline{a}\tau)$.
\end{corol}

\pr By a result in \cite[Corollary 3.3.1]{S} the function
$f(z)=z-F(z), \, z\in\Delta$  is a generator of a one-parameter
semigroup. By our assumptions we have $f(\tau)=0$ and $f(0)=-a$.
Hence, by Corollary 2 $\, f'(\tau)\leq -2\Re (\overline{a}\tau)$,
and $F'(\tau)\geq 1+2 \Re (\overline{a}\tau)$. \epr

Now let us consider a class of functions
$f\in\mathrm{Hol}(\Delta,\mathbb{C})$ which are continuous on
$\overline{\Delta}$ and satisfy the boundary condition
\begin{equation}\label{a18}
\Re f(z)\overline{z}\geq |f(z)|\cos\frac{\alpha\pi}{2} \quad
\mbox{for all}\quad z\in\partial\Delta,
\end{equation}
for some $\alpha\in (0,2]$. As we already mentioned if $\alpha\leq
1$ then condition (\ref{a18}) implies $f\in\G (\Delta)$ (cf.
Proposition 1 (ii)). Conversely, if $f\in\G(\Delta)$ is continuous
on $\overline\Delta$, then (\ref{a18}) holds with $\alpha=1$. So,
this class generalize in a sense  the class of holomorphic
generators which are continuous on  $\overline\Delta$.

\begin{pred0}
Let $f\in\mathrm{Hol}(\Delta,\mathbb{C})$ be continuous on
$\overline{\Delta}$  and satisfy the condition (\ref{a18}). Then
the condition
\begin{equation}\label{a19}
\lim\limits_{\begin{smallmatrix}z\to\tau \\
z\in\overline\Delta\end{smallmatrix}}
\frac{f(z)}{(z-\tau)^{2+\alpha}}=0 \quad\mbox{ for some} \quad
\tau\in\partial\Delta
\end{equation}
implies that $f\equiv 0$.
\end{pred0}

\pr Denote $$g(z)=\frac{f(z)}{(z-\tau)(1-\overline{\tau} z)} \,
.$$ The continuity of $f$ and (\ref{a19}) imply that this function
is continuous (consequently, bounded) on $\overline{\Delta}$.

Now we rewrite (\ref{a18}) in the form: $$-\Re \left[
\overline{\tau}(\tau-z)^{2}g(z)\overline{z} \right ]\geq
|\tau-z|^{2}\cdot |g(z)|\cdot \cos\frac{\alpha\pi}{2}, \quad
z\in\partial\Delta.$$  Hence, $$\Re g(z)\geq |g(z)|\cdot \cos
\frac{\alpha\pi}{2}, \quad z\in\partial\Delta\setminus\{\tau\}.$$
This inequality also holds at the point $\tau$ because of the
continuity of $g$.

It follows from the subordination principle for subharmonic
functions (see, for example, \cite[p. 396]{G}) that the latter
inequality holds for all $z\in\overline{\Delta}$. Geometrically
this fact means that $g$ maps $\Delta$ into the sector
$\overline{A_{\alpha}}$, where $$A_{\alpha}=\left \{
w\in\mathbb{C}:|\arg w|<\frac{\alpha\pi}{2}, \; \alpha\in(0,2]
\right \}.$$

Suppose that there exists $z\in\Delta$ such that
$w=g(z)\in\partial A_{\alpha}$. Then by the maximum principle
$g\equiv{\rm const} =w$ and $f(z)=w\tau(z-\tau)^{2}$. In this case
$w$ must be zero, since otherwise we get contradiction
with~(\ref{a19}). Hence, either $w=0$ or $g(\Delta)\subset
A_{\alpha}$.

If $w=0$ then $f\equiv 0$ and we are done.

Let now $g(\Delta)\subset A_{\alpha}$. Equality (\ref{a19})
implies that
$$\angle\lim\limits_{z\rightarrow
\tau}\frac{g(z)}{(z-\tau)^{\alpha}}=-\tau\angle\lim\limits_{z\rightarrow
\tau}\frac{f(z)}{(z-\tau)^{2+\alpha}}=0.$$

Applying Theorem C we get $g\equiv 0$, hence $f\equiv 0$ . \epr

\begin{corol}
Let $F\in\mathrm{Hol}(\Delta,\mathbb{C})$ be continuous on
$\overline{\Delta}$  and satisfy the boundary condition
\begin{equation}\label{a20}
\Re F(z)\overline{z}\leq 1-|F(z)-z|\cos\frac{\alpha\pi}{2},\quad
z\in\partial\Delta,
\end{equation}
for some $\alpha\in(0,2]$. If there exists $\tau\in\partial\Delta$
such that
$$F(z)=\tau+(z-\tau)+o(|z-\tau|^{2+\alpha})$$ when $z\rightarrow\tau$, then $F\equiv I.$
\end{corol}

\section{Commuting semigroups.}

\begin{pred0}
Let $f$ and $g$ be generators of one-parameter commuting
semigroups $\left\{ F_{t} \right\}_{t \geq 0}$ and $\left\{ G_{t}
\right\}_{t \geq 0}$, respectively, and $f(\tau)=0$ at some point
$\tau\in\overline{\Delta}$.

(i) Let $\tau\in\Delta$. If $f'(\tau)=g'(\tau)$ then $f\equiv g$.

(ii) Let $\tau\in\partial\Delta$. Suppose $f$ and $g$ admit the
following representations
\begin{equation}\label{star1}
f(z)=f'(\tau)(z-\tau)+\ldots +\frac{f^{(m)}(\tau)}{m!}(z-\tau)^m
+o\left(|z-\tau|^m\right)
\end{equation}
and
\begin{equation}\label{star2}
g(z)=g(\tau)+g'(\tau)(z-\tau)+\ldots
+\frac{g^{(m)}(\tau)}{m!}(z-\tau)^m +o\left(|z-\tau|^m\right)
\end{equation}
when $z\to\tau$ along some curve lying in $\Delta$ and ending at
$\tau$. If $f^{(m)}(\tau)=g^{(m)}(\tau)\neq 0$, then $f\equiv g$.
\end{pred0}

\begin{rem}
If $\tau\in\partial\Delta$ is the Denjoy--Wolff point of a
semigroup generated by a mapping $h\in\G(\Delta)$, then $h$ admits
the expansion
\[
h(z)=h'(\tau)(z-\tau)+o(z-\tau)
\]
when $z\to\tau$ in each non-tangential approach region at $\tau$
and $h'(\tau)=\angle\lim\limits_{z\to\tau}h'(z)$. Moreover, in
this case $h'(\tau)$ is a non-negative real number which is zero
if and only if $h$ generates a semigroup of parabolic type (see
\cite{E-S}).

Therefore, if $f$ (or $g$) in Theorem~3 generates a semigroup of
hyperbolic type with the Denjoy--Wolff point
$\tau\in\partial\Delta$ then the condition $f'(\tau)=g'(\tau)$ is
enough to provide that $f\equiv g$.
\end{rem}

\begin{rem}
As a matter of fact, if $f$ and $g$ have expansion (\ref{star1})
and (\ref{star2}) when $z\to\tau$ in each non-tangential approach
region at $\tau\in\partial\Delta$ up to the third order $m=3$,
such that $f'(\tau)=g'(\tau),\ f''(\tau)=g''(\tau)$ and
$f'''(\tau)=g'''(\tau)$ then $f\equiv g$.

If, in particular, $f^{(i)}(\tau)=g^{(i)}(\tau)=0,\ i=1,2,3,$ then
both $f$ and $g$ are equal zero identically by Theorem~1.
\end{rem}

Theorem~3 is a consequence of the following more general
assertion.

Define two linear semigroups $\left\{ A_{t} \right\}_{t \geq 0}$
and $\left\{ B_{t} \right\}_{t \geq 0}$  of composition operators
on $\mathrm{Hol}(\Delta,\mathbb{C})$ by
\begin{equation}\label{a39}
A_{t}(h)=h\circ F_{t}\quad\mbox{and}\quad B_{t}(h)=h\circ G_{t},
\quad t\geq 0.
\end{equation}
The operators $\Gamma_f$ and $\Gamma_g$ defined by
\begin{equation}\label{a340}
\Gamma_{f}(h)=h'f\quad\mbox{and}\quad \Gamma_{g}(h)=h'g
\end{equation}
are their generators, respectively.

\begin{pred0}
Let $f$ and $g\in\mathrm{Hol}(\Delta,\mathbb{C})$ be generators of
one-parameter  semigroups $\left\{ F_{t} \right\}_{t \geq 0}$ and
$\left\{ G_{t} \right\}_{t \geq 0}$,  respectively. Let $A_{t}$
and $B_{t}$ be defined by (\ref{a39}). Then the following are
equivalent:

(i) $F_{t}\circ G_{s}=G_{s}\circ F_{t}$, $s,t\geq 0$, i.e., the
semigroups $\left\{ F_{t} \right\}_{t \geq 0}$ and $\left\{ G_{t}
\right\}_{t \geq 0}$ are commuting;

(ii) $A_{t}\circ B_{s}=B_{s}\circ A_{t},\ s,t\geq 0$, i.e., the
linear semigroups $\left\{ A_{t} \right\}_{t \geq 0}$ and $\left\{
B_{t} \right\}_{t \geq 0}$ are commuting;

(iii) $\Gamma_{f}\circ \Gamma_{g}=\Gamma_{g}\circ \Gamma_{f} $,
i.e., the linear semigroup generators $\Gamma_{f}$ and
$\Gamma_{g}$ are commuting;

(iv) the Lie commutator $$[f,g]=f'g-g'f=0;$$

(v) $f=\alpha g$ for some $\alpha\in\mathbb{C}$.
\end{pred0}

\pr Suppose that $f\not\equiv 0$. First we prove the equivalence
of assertions (i) and (v).

Let (i) holds. If $f(\tau)=0$, $\tau\in\Delta$, then $\tau$ is a
unique common fixed point for the semigroup $\left\{ F_{t}
\right\}_{t \geq 0}$ generated by $f$, i.e., $F_{t}(\tau)=\tau$
for all $t \geq 0$ (see, for example, \cite{B-P}, \cite{S} ).

If $F_{t}$ and $G_{s}$ are commuting for all $s,t \geq 0$, then we
have $$G_{s}(\tau)=G_{s}\left (F_{t} (\tau) \right )=F_{t}\left
(G_{s} (\tau) \right ).$$ Hence, it follows by the uniqueness of
the fixed point $\tau$ that $G_{s} (\tau)=\tau$ for all $s\geq 0$,
and so $g(\tau)=0$.

Consider the function $h\in\mathrm{Hol}(\Delta,\mathbb{C})$
defined by the differential equation
\begin{equation}\label{a22}
 \mu h (z) =h'(z)f(z).
\end{equation}
It is known that if $\mu=f'(\tau) $ then equation (\ref{a22}) has
a unique solution $h\in\mathrm{Hol}(\Delta,\mathbb{C})$ normalized
by the condition $h'(\tau)=1$ (see \cite{S}).

In addition, this function $h$ solves Schroeder's functional
equation
\begin{equation}\label{a23}
 h\left ( F_{t}(z)\right )=e^{-\mu t}h(z).
\end{equation}
Now, for any $s,t\geq 0$ we get from (\ref{a23})  $$h\left
(G_{s}\left (F_{t} (z) \right )\right )=h\left (F_{t}\left (G_{s}
(z) \right )\right )=e^{-\mu t}h\left (G_{s}(z)\right ).$$ Denote
$h_{s}=h\circ G_{s}$. Then we have
\begin{equation}\label{a24}
 h_{s}\left ( F_{t}(z)\right )=e^{-\mu t}h_{s}(z).
\end{equation}
Differentiating (\ref{a24}) at $t=0^{+}$ we get
\begin{equation}\label{a55}
\mu h_{s}(z)=h'_{s}(z)f(z).
\end{equation}
Comparing (\ref{a22}) and (\ref{a55}) implies $h_{s}(z)=\lambda
(s)h(z)$ for some $\lambda (s)\in \mathbb{C}$, or
\begin{equation}\label{a25}
 h\left ( G_{s}(z)\right )=\lambda (s)h(z).
\end{equation}

Since the left-hand side of the latter equality is differentiable
in $s\geq 0$, the scalar function $\lambda (s)$ is differentiable
too. Differentiating (\ref{a25}) at $s=0^{+}$ we get
\begin{equation}\label{a26}
 \lambda '(0)h(z)=-h'(z)g(z).
\end{equation}
Note that $h(\tau)=0$ while $h(z)\neq 0$ for all $z\in \Delta,
z\neq\tau$. In addition, it can be shown (see \cite{S}) that $h$
is univalent. Hence, $h'(z)\neq 0$ for all $z\in\Delta$.

Finally, we obtain from (\ref{a22}) and (\ref{a26}) that
\[
f(z)=\alpha g(z), \quad\mbox{where} \quad
\alpha=-\frac{\mu}{\lambda '(0)}\,.
\]

Now, let us suppose that $f$ has no null point in $\Delta$. Then
the function $p:\Delta\mapsto\mathbb{C}$ given by
\begin{equation}\label{a27}
p(z)=-\int\limits_{0}^{z}\frac{d\varsigma}{f(\varsigma)}
\end{equation}
is well defined holomorphic function on $\Delta$ with $p(0)=0$.

Recall that the semigroup $\left\{ F_{t} \right\}_{t \geq 0}$
generated by $f$ can be defined by the Cauchy problem
\begin{equation}\label{a28}
\left \{ \begin{smallmatrix}\dfrac{dF_{t}(z)}{dt} \displaystyle
+f\left ( F_{t}(z)\right )=0, \, t\geq 0\\
\!\!\!\!\!\!\!\!\!\!\!\!\!\!\!\!\!\!\!\!\!\!\!\!\!\!\!\!\!\!\!\!\!\!\!\!\!\!\displaystyle
F_{0}(z)=z, \quad z\in\Delta
\end{smallmatrix}\right.
\end{equation}
Substituting here $f(z)=-\frac{1}{p'(z)}$ we obtain $$p'\left (
F_{t}(z)\right ) dF_{t}(z)=dt.$$ Integrating the latter equality
on the interval $[0,t]$ we get that $p$ is a solution of Abel's
functional equation
\begin{equation}\label{a29}
 p\left ( F_{t}(z)\right )=p(z)+t.
\end{equation}
Now, for any fixed $s\geq 0$ we have $$p\left (G_{s}\left (F_{t}
(z) \right )\right )=p\left (F_{t}\left (G_{s} (z) \right )\right
)=p\left (G_{s} (z) \right )+t.$$ Once again, setting
$p_{s}=p\circ G_{s}$, we have
\begin{equation}\label{a30}
 p_{s}\left ( F_{t}(z)\right )=p_{s}(z)+t.
\end{equation}
Differentiating (\ref{a30}) at $t= 0^{+}$ we get
\begin{equation}\label{a31}
 p'_{s}(z)=-\frac{1}{f(z)}\,,
\end{equation}
and by (\ref{a27}), \, $p_{s}(z)=p(z)+\kappa (s), \; \kappa
(s)\in\mathbb{C}$, or
\begin{equation}\label{a32}
p\left (G_{s}(z)\right )=p(z)+\kappa (s).
\end{equation}
Differentiating (\ref{a32}) at $s= 0^{+}$ we obtain the equality
\begin{equation}\label{a42}
p'(z)=-\frac{\kappa '(0)}{g(z)}.
\end{equation}
Comparing (\ref{a31}) and (\ref{a42}) gives
\begin{equation}\label{a33}
f=\alpha g \quad\mbox{with} \quad \alpha=\frac{1}{\kappa '(0)}\,.
\end{equation}

Now we prove that (v)$\Rightarrow$(i). Let $f=\alpha g$ for some
$\alpha\in\mathbb{C}$.

First we assume that $g$ has an interior null-point
$\tau\in\Delta$. In this case there is a univalent solution of the
differential equation
\begin{equation}\label{a34}
\mu h(z)=h'(z)g(z)
\end{equation}
with some $\mu\in\mathbb{C}, \quad \Re\mu\geq 0$.

Since $f=\alpha g$, we have that $h$ is also a solution of the
equation
\begin{equation}\label{a35}
\nu h(z)=h'(z)f(z), \quad \nu =\alpha\mu.
\end{equation}

In turn, equations (\ref{a34}) and (\ref{a35}) are equivalent to
Schroeder's functional equations
\begin{equation}\label{a36}
h\left (G_{s} (z) \right )=e^{-\mu s}h(z), \quad s\geq 0
\end{equation}
and
\begin{equation}\label{a37}
h\left (F_{t} (z) \right )=e^{-\nu t}h(z), \quad t\geq 0, \quad
\nu = \alpha\mu,
\end{equation}
respectively, where $\left\{ F_{t} \right\}_{t \geq 0}$ is the
semigroup generated by $f$.

Consequently,
\begin{eqnarray*}
F_{t}\left (G_{s} (z) \right )=h^{-1}\left (e^{-\nu t}h\left
(G_{s}(z)\right )\right )=h^{-1}\left (e^{-\nu t}\cdot e^{-\mu
s}h(z)\right )\\ =h^{-1}\left ( e^{-\mu s}h\left (F_{t}(z)\right
)\right )=G_{s}\left (F_{t} (z) \right )
\end{eqnarray*}
for all $s,t\geq 0$ and we are done.

Now let us assume that $g$ has a boundary null-point
$\tau\in\partial\Delta$ with $g'(\tau)\ge0$ (see Remark 1 above).
In this case for each $c\in\mathbb{C}, \, c\neq 0$, Abel's
equations
\[
p\left (G_{s}(z)\right )=p(z)+c s
\]
and
\[
p\left (F_{t}(z)\right )=p(z)+c\alpha t
\]
have the same solution
\[
p(z)=-c\int\limits_{0}^{z}
\frac{d\varsigma}{g(\varsigma)}=-c\alpha\int\limits_{0}^{z}
\frac{d\varsigma}{f(\varsigma)},
\]
which is univalent on $\Delta$.

Once again we calculate
\begin{eqnarray*}
F_{t}\left (G_{s}(z)\right )=p^{-1}\left (p\left (G_{s}(z)\right
)+c \alpha t \right )=p^{-1}\left (p(z)+c \alpha t + c s\right )\\
=p^{-1}\left (p\left (F_{t}(z)\right )+ c s \right )=G_{s}\left
(F_{t}(z) \right ).
\end{eqnarray*}
The implication (v)$\Rightarrow$(i) is proved.

The equivalence of (i) and (ii) is obvious.

To verify the equivalence of (iii) and (iv) we just calculate:
$$\Gamma_{f}\left (\Gamma_{g}(h) \right )=h''g f+h'g'f,$$
 $$\Gamma_{g} \left ( \Gamma_{f} (h) \right ) = h''f g+h'f'g.$$
Hence, $\Gamma_{f}\circ \Gamma_{g}=\Gamma_{g}\circ \Gamma_{f} $ if
and only if $f'g-g'f=0$.

Now, it is clear, that (v) implies (iv).

Finally we prove the implication (iv) $\Rightarrow$ (v).
Obviously, (iv) implies that if $f$ has no null points in $\Delta$
then $g$ also has no null points in $\Delta$ and, hence, (v)
follows. If $f(\tau)=0$ for some $\tau\in\Delta$, then also
$g(\tau)=0$, and by (\ref{a43}) one can write $f(z)=(z-\tau)p(z)$
and  $g(z)=(z-\tau)q(z)$, where $p$ and $q$ do not vanish in
$\Delta$. Now it follows that $$[f,g]=(z-\tau)[p,q]$$  Hence,
again we have $p=aq$, and hence $f=ag$ for some $a\in\mathbb{C},
\, a\neq 0$. \epr

\prrr First we note, that by Theorem 4
\begin{equation}\label{a38}
f=\alpha g, \quad \alpha\in\mathbb{C}.
\end{equation}

(i) Let $f'(\tau)=g'(\tau)=0$. By Proposition 1 $\; f$ admits
representation  $$f(z)=(z-\tau )(1-\overline{\tau}z)p(z),\;\;
z\in\Delta,$$ where $p\in\mathrm{Hol}(\Delta,\mathbb{C})$, $\Re
p(z)\geq 0$.

Since $f'(\tau)=(1-|\tau|^{2})p(\tau)=0$, we have $p(\tau)=0$ and
it follows from the maximum principle that $p\equiv 0$. Hence,
$f\equiv 0$ and by (\ref{a38}) also $g\equiv 0$.

Assume now $f'(\tau)=g'(\tau)\neq 0$. Then it follows from
(\ref{a38} ) that $\alpha=1$ and so $f\equiv g$.

(ii) In general, by (\ref{a38}) we have $f^{(k)}(\tau)=\alpha
g^{(k)}(\tau)$, $0<k\leq m$. Hence, the condition
$f^{(k)}(\tau)=g^{(k)}(\tau)\neq 0$ for some $0<k\leq m$ implies
that $\alpha=1$ and, consequently, $f\equiv g$. \epr

Let $\mathcal{S}_f=\left\{F_t\right\}_{t\ge0}$ be the semigroup
generated by $f\in\G(\Delta)$. The set
$\mathcal{Z}(\mathcal{S}_f)$ of all semigroups
$\mathcal{S}=\left\{G_t\right\}_{t\ge0}$ such that
\[
F_t\circ G_s=G_s\circ F_t,\quad t,s\ge0,
\]
is called the {\bf centralizer of $\mathcal{S}_f$}.

It is clear that for each $f\in\G(\Delta)$ the centralizer
$\mathcal{Z}(\mathcal{S}_f)$  contains  $\mathcal{S}_{\alpha f}$
for all $\alpha\ge0$.

Therefore we will say that {\bf the centralizer of $\mathcal{S}_f$
is trivial} when the conclusion
$\mathcal{S}\in\mathcal{Z}(\mathcal{S}_f)$ implies that
$\mathcal{S}=\mathcal{S}_{\alpha f}$ for some $\alpha\ge0$.

\begin{propo}
Let $f$ be the generator of a semigroup
$\mathcal{S}_f=\left\{F_t\right\}_{t\ge0}$, and let
$\tau\in\partial\Delta$ be the Denjoy--Wolff point of
$\mathcal{S}_f$. Then if one of the following conditions holds
then the centralizer $\mathcal{Z}(\mathcal{S}_f)$ is trivial:

(i) $\mathcal{S}_f$ is a hyperbolic type semigroup ($f'(\tau)>0$)
which is not a group;

(ii) $f$ admits the expansion
\[
f(z)=a(z-\tau)^3+o\left((z-\tau)^3\right)\quad \mbox{with}\quad
a\not=0
\]
when $z\to\tau$ in each non-tangential approach region at $\tau$.
\end{propo}

The first statement is based on the following simple lemma.

\begin{lemma}
Let $f$ and $g$ be generators of two nontrivial (neither $f$ nor
$g$ are identically zero) commuting semigroups
$\mathcal{S}_f=\left\{F_t\right\}_{t\ge0}$ and
$\mathcal{S}_g=\left\{G_t\right\}_{t\ge0}$, respectively. Then
$\mathcal{S}_f$ is of hyperbolic type if and only if
$\mathcal{S}_g$ is. In this case $f=\alpha g$ with real $\alpha$.
Moreover, $\alpha<0$ implies that $\mathcal{S}_f$ and
$\mathcal{S}_g$ are both groups of hyperbolic automorphisms having
`opposite' fixed points, i.e., the attractive point for
$\mathcal{S}_f$ is the repelling point for $\mathcal{S}_g$ and
conversely.
\end{lemma}

\pr Since  $\mathcal{S}_f$ and $\mathcal{S}_g$ are commuting then
by Theorem~4 there exists $\alpha\in\mathbb{C}$ such that
$f=\alpha g$. In our settings $\alpha$ is not zero.  If $\tau$ is
the Denjoy--Wolff point of $\mathcal{S}_f$ then $f(\tau)=0$ and
therefore also $g(\tau)=0$. Now since  $f'(\tau)>0$ then
$g'(\tau)=\frac1\alpha\,f'(\tau)$ exists finitely and it must be a
real number by Corollary~2. So must be $\alpha$.

Now let us assume that $\alpha$ is negative. Then
$g'(\tau)=\frac1\alpha\,f'(\tau)<0$. Hence the semigroup
$\mathcal{S}_g$ generated by $g$ must have the Denjoy--Wolff point
$\sigma\in\overline\Delta$ different from $\tau$.

It is clear that $\sigma$ can not be inside $\Delta$ since
otherwise it must be a common fixed point of both semigroups
$\mathcal{S}_f$ and $\mathcal{S}_g$ because of the commuting
property.

So, $\sigma\in\partial\Delta$ and $g'(\sigma)\ge0$ (see
\cite{E-S}), then $f(\sigma)=0$ and $f'(\sigma)\le0$. It follows
by a result in \cite{SD1} that
\begin{equation}\label{repr}
0<f'(\tau)\le-f'(\sigma)
\end{equation}
and the equality is possible if and only if $f$ is the generator
of a group of hyperbolic automorphisms. By the same theorem we
have the reversed inequality for $g$
\[
0\le g'(\sigma)\le- g'(\tau)
\]
that means
\[
0\le \frac1\alpha\, f'(\sigma)\le-\frac1\alpha\, f'(\tau).
\]
Comparing this inequality with (\ref{repr}) gives us that
$f'(\tau)=-f'(\sigma)>0$ and $g'(\tau)=-g'(\sigma)<0$ which means
that both $f$ and $g$ generate groups of hyperbolic automorphisms
with opposite fixed points. \epr

\noindent{\bf Remark.} The last assertion of this lemma follows
also by a result of Behan (see \cite{Be}). Indeed, let $\alpha<0$.
Then the equality $f(z)=\alpha g(z)$ implies that $g'(\tau)$
exists and is a real negative number. So, the Denjoy--Wolff point
$\tau$ of the semigroup $S_f$ can not be the Denjoy--Wollf point
of the semigroup $S_g$. Hence by \cite{Be} we conclude that $S_f$
and $S_g$ are groups of hyperbolic automorphisms.

\noindent{\bf Proof of Proposition 1.} The statement (i) is a
direct consequence of the previous lemma. To prove the second
statement we note that by Theorem~1 the number $a\tau^2$ is a
non-negative real number. On the other hand, since $\mathcal{S}_f$
and $\mathcal{S}_g$ commute by Theorem~4 there is a number
$\alpha\in\mathbb{C}$ such that $f=\alpha g$.

Therefore, since $\alpha\not=0$ also $g$ admits the expansion
\[
g(z)=\frac a\alpha\,(z-\tau)^3+o\left((z-\tau)^3\right)
\]
and again by Theorem~1 we have that also $\frac
a\alpha\,\tau^2\ge0$. This implies that $\alpha$ is a nonnegative
real number. \epr

A natural question which arises in the context of the above
theorem is:

\noindent$\bullet$ If two elements $F_{t_0}$ and $G_{s_0}$ of
semigroups $S_{f}=\{F_{t}\}_{t\geq 0}$ and $S_{g}=\{G_{t}\}_{t\geq
0}$ commute for some positive $t_{0}$ and $s_{0}$, do these
semigroup $S_{f}$ and $S_{g}$ commute in the sense:
\[
F_{t}\circ G_{s}=G_{s}\circ F_{t}
\]
for each pair $t,s\geq 0.$

The answer is immediately affirmative due to a more general result
of C.~C.~Cowen (\cite{CC}), Corollary \ ) if neither $F_{t_{0}}$
nor $G_{s_{0}}$, respectively, are of parabolic type.

The situation becomes more complicated if $F_{t_{0}}$,
respectively $G_{s_{0}}$, are parabolic.

Example 4.4 in \cite{CC} shows that there is a triple of such
mappings $F,\ G_{1}$ and $G_{2}$ in $\Hol(\Delta)$ for which
$G_{1}$ and $G_{2}$ commute with $F$, but they do not commute each
other.

Nevertheless, under some additional requirements on smoothness at
the boundary Denjoy-Wolff point repeating the arguments using in
the proof of Theorem 1.2 in \cite{T} one can give an affirmative
answer the above question. Namely,

\noindent$\bullet$ Let $F_{t_{0}}$ and $G_{s_{0}}$ be two
commuting elements of semigroups $S_{f}$ and $S_{g}$,
respectively, $t_{0},s_{0}>0$, and let $F_{t_{0}}$ is of parabolic
type with a Denjoy--Wolff point $\tau\in\partial\Delta.$ If both
$F_{t_{0}}$ and $G_{s_{0}}$ belong to the class $C^{2}(\tau )$ and
$F''_{t_{0}}(\tau)$ as well as $G''_{t_{0}}(\tau)$ do not vanish,
then $f=ag$ for some $a\in\mathbb{C}$, i.e., the semigroups
$S_{f}$ and $S_{g}$ commute:
\[
F_{t}\circ G_{s}=G_{s}\circ F_{t}
\]
for all \ $t,s\geq 0.$

\end{document}